%% file: articlefinalnew.tex
\documentclass[a4paper,12pt]{article}
\usepackage[toc,page]{appendix}
\usepackage{fancyhdr,lastpage}
\usepackage{amsthm}
\usepackage{amsmath}
\usepackage{amsfonts}
\usepackage{fullpage}
\usepackage{amssymb}
\numberwithin{equation}{section}
\newtheorem{theor}{Theorem}[section]
\newtheorem{propo}{Proposition}[section]
\newtheorem{lem}{Lemma}[section]
\newtheorem{definition}{Definition}[section]

\newtheorem{rem}{Remark}[section]
\usepackage{bbm}
\numberwithin{equation}{section}
\usepackage{graphics}
\usepackage{url}
\usepackage{enumerate}
\usepackage{hyperref}
\usepackage{pgf}

\begin{document}
\date{}
\title{Well-posedness and gradient blow-up estimate near the boundary
for a Hamilton-Jacobi equation with degenerate diffusion}
\author{Amal Attouchi}
\maketitle

\begin{abstract}
This paper is  concerned with weak solutions of the degenerate viscous Hamilton-Jacobi equation
$$\partial_t u-\Delta_p u=|\nabla u|^q,$$ with Dirichlet boundary conditions in a bounded domain $\Omega\subset\mathbb{R}^N$, where $p>2$ and $q>p-1$. With the goal of studying the gradient blow-up phenomenon for this problem,  we first establish  local well-posedness with blow-up alternative  in $W^{1, \infty}$ norm. We then obtain a precise gradient estimate involving the distance to the boundary. It shows in particular that the gradient blow-up can take place only on the boundary. A regularizing effect for $u_t$ is also obtained.
\end{abstract}

\section{Introduction and main results}

  This article is concerned with the existence and qualitative properties of  weak solutions of the initial boundary value problem of  the $p$-Laplacian with  a nonlinear gradient  source term
\begin{equation}\label{equi}
 \left\{
\begin{array}{lll}
 \partial_t u- \mathrm{div}( |\nabla u|^{p-2} \nabla u)=|\nabla u|^q,  &  x\in\Omega, t>0,\\
u(x,t)=g(x),                                                  & x\in\partial\Omega, t>0,\\
u(x,0)=u_0(x),                                             & x\in\Omega,

\end{array}
\right.
\end{equation}
 where $\Omega$ is a bounded domain in $\mathbb{R}^N$ of class $C^{2+\alpha}$ for some $\alpha>0$, $p>2$ and $q>p-1$.
Throughout the paper we assume that  the boundary data $g\geq 0$ is the trace on $ \partial\Omega$  of a regular function  in $ C^2(\overline{\Omega})$,  also denoted $g$,  and  the initial data  $u_0$  satisfies
\begin{equation}\label{initial}
u_0\in W^{1,\infty}(\Omega),\qquad  u_0\geq 0,\quad u_0(x)=g(x) \qquad\text{for}\quad x\in\partial\Omega.
\end{equation}
We note that, as far as bounded solutions are concerned, there is no loss of generality in assuming $g, u_0\geq 0$, since the partial differential equation in \eqref{equi}  is unchanged  when adding a constant to $u$.

When $p=2$, the differential equation of \eqref{equi} is the so-called viscous Hamilton-Jacobi equation and it appears in the physical theory of growth and roughening of surfaces, where it is known as the Kardar-Parisi-Zhang equation ($q=2$), and has been studied by many authors (see for example \cite{local, superlinear} and the references therein). It is known  that, under certain
conditions, $|\nabla u|$ blows up in a finite time $t = T_{max}$ while, by the maximum principle, all
solutions are uniformly bounded (cf. \cite{souma, hesa, Zhan}). We shall call such phenomenon gradient blow-up (GBU). This is different from the usual blow-up in which the $L^{\infty}$ norm of the solution tends to infinity
as $t\to T_{max}$ (cf. \cite{superlinear}).
Sharp results on  gradient blow-up analysis, including blow-up rate, blow-up set, blow-up profile and continuation after blow-up have been recently obtained, see e.g. \cite{single, guo,hesa, superlinear,Arri,vaz} and the references therein.

When $p>2$, equation \eqref{equi} is a degenerate parabolic equation for $| \nabla u|=0$ and one cannot expect the existence of classical solutions. Weak solutions can be obtained by approximation with solutions of regularized problems. This was done in \cite{zhao2} when the right hand side in \eqref{equi} is replaced with a general nonlinearity $f(u,\nabla u,x,t)$. In the case where $f$ depends on $\nabla u$, typically for problem \eqref{equi}, the results in \cite{zhao2} require the assumption $q\leq p-1$,
in which case a global solution is directly constructed for any initial data.
Local-in-time existence results are also given in \cite{zhao2}  but they require that $f$ actually does not depend on $\nabla u$. In \cite{chen2}, the existence of a global weak solution for $q>p-1$ was proved for small data,
under the assumption that the mean curvature of $\partial\Omega$ is nonpositive. In the articles \cite{oups, barle}, problem \eqref{equi} was studied in the framework of viscosity solutions, but only in situations where global existence of a $W^{1, \infty}$ solution is guaranteed, namely for $q\leq p$ or for suitably small initial data when $q>p$.  On the other hand,  when $q>p$,  global existence  is not  expected in general for large initial data. A result in this direction was given in [\cite{oups} Theorem 5.2], where it was proved that problem \eqref{equi} (with $g=0$)  cannot admit a global, Lipschitz continuous, weak solution for large initial data. See \cite{lieb, dlot, giga} and the references therein for earlier counter-examples concerning related quasilinear equations.

Our first goal will be to complete the above results by constructing a unique, maximal in time, $W^{1, \infty}$ solution, without size restriction on the initial data and to establish the blow up alternative in $W^{1,\infty}$ norm. This will enable us to interpret the above mentioned  global nonexistence result from \cite{oups} appropriately as a gradient blow-up (GBU) result (see Theorem\ref{blo} and Remark \ref{secblo} below), and will provide  the grounds for the subsequent analysis of the asymptotic behavior of GBU solutions. For the local existence part, we will follow and suitably modify the approximation procedure used in \cite{zhao2}.

The main difficulty is to get relevant estimates on the first order derivatives of the approximate solutions in order to pass to the limit in the nonlinear source term. To deal  with this difficulty, our main new ingredient with respect to \cite{zhao2} is the construction of suitable barrier functions,  in order to get uniform pointwise estimates on the gradients near the boundary for small time.
We then use a strong result  of DiBenedetto  and Friedman \cite{dib} on the $\mathrm{H\ddot{o}}$lder regularity of  gradients of  weak solutions of degenerate parabolic equations and consequently we will use the framework of weak rather than viscosity solutions.

 First, let us state the precise definition of solution. Let $Q_T=\Omega\times (0,T)$ and $\partial_p Q_T=\left\{\partial\Omega\times[0, T]\right\}\cup\left\{\overline{\Omega}\times\left\{0\right\}\right\}$, $T>0$. Throughout this paper, we will use the following definition of weak solution for \eqref{equi}.
\begin{definition}\label{def1}
Set $m=max(p,q)$.  A function  $u(x,t)$ is called a weak super- (sub-) solution of problem \eqref{equi} on $Q_T$ if
\begin{eqnarray*}
 &u\in C(\overline{\Omega}\times [0,T))\cap L^m((0,T); W^{1,m}(\Omega)),\\
& u_t\in L^2((0,T); L^2(\Omega)),\\
&u(x,0)\geq (\leq)\, u_0(x),\qquad   u\geq  (\leq)\, g \,\text{on}\,\, \partial\Omega \, \text{and}
\end{eqnarray*}
\begin{equation}
 \int\int_{Q_T} u_t \psi+|\nabla u|^{p-2} \nabla u\cdot\nabla \psi \, dx\, dt\geq(\leq)\int\int _{Q_T} |\nabla u|^q \psi \,dx \,dt
\end{equation}
holds for all $\psi\in C^0(\overline{Q_T}) \cap L^p((0,T); W^{1,p}(\Omega))$ such that $ \psi\geq 0$, $\psi=0$ on $\partial\Omega\times(0,T)$.
A function $u$ is a weak solution of \eqref{equi} if it is a super-solution and a sub-solution.
\end{definition}

Our first result concerns local existence and uniqueness of  weak solutions (see also Section 2  for a comparison principle).

\begin{theor}\label{exis}
Assume that $q> p-1>1$. Let $M>0$ and let  $u_0$ satisfy \eqref{initial} and $\left\|\nabla u_0\right\|_{\infty}\leq M$. Then 

\begin{enumerate}[(i)]
\item  There exist a time $T=T\left(M, p, q, N, \left\|g\right\|_{C^2}\right)>0$
 and a weak solution  $u$  of (\ref{equi}) on $[0,T)$, which moreover satisfies $u\in L^\infty_{loc} ([0,T);W^{1,\infty}(\Omega))$.
 
\item For any $\mathcal{T}>0$ the problem  (\ref{equi}) has at most one weak solution $u$ such that
$u\in L^\infty_{loc}([0,\mathcal{T});W^{1,\infty}(\Omega))$.

\item  
There exists a (unique) maximal, weak solution of \eqref{equi}, still denoted by $u$. Let  $T_{max}(u_0)$ be its  existence time.

Then  
\begin{equation}\label{pmax}
\underset{\Omega}{\text{min}}\; u_0\leq u \leq\underset{\Omega}{\text{max}}\; u_0 \qquad \text{in}\quad \Omega\times\left( 0, T_{max}(u_0)\right)
\end{equation}
 and  $$\text{if}\quad T_{max}(u_0)<\infty,\quad\text{then}\quad \underset{t\to T_{max}(u_0)}{\mathrm{lim}} \ \left\|\nabla u\right\|_{L^{\infty}(\Omega)}=\infty\qquad  \text{(gradient blow up GBU)}.$$
\end{enumerate}
\end{theor}
\begin{rem}\label{rem1}
Concerning Definition \ref{def1}, we note that if $0<T_1<T_2<\infty$ and $u$ is a weak solution on $Q_{T_2}$, then the restriction of $u$ to $Q_{T_1}$ is a weak solution on $Q_{T_1}$ (this can be easily checked, taking any test function $\psi$ on $Q_{T_1}$, by extending $\psi$ as $\tilde{\psi}_n(x,t)=\psi(x,T_1)[1-n(t-T_1)]_+$ for $t\in (T_1,T_2]$ and letting $n\to\infty$).
Then, in Theorem \ref{exis}(iii), by $u$ being the maximal weak solution of \eqref{equi}, we mean that
$u$ is a weak solution on $Q_\tau$ for any $\tau\in (0,T_{max}(u_0))$ but cannot be extended to a weak solution on $Q_{T'}$ for any $T'>T_{max}(u_0)$.
\end{rem}
We next establish a precise gradient estimate involving the distance to the boundary. Here and in the rest of the paper we denote $\delta(x)=\text{dist} (x, \partial\Omega)$.
\begin{theor}\label{profil}
Let  $q> p-1>1$. Let $M>0$ and let $u_0$ satisfy \eqref{initial}  and  $\left\|\nabla u_0\right\|_{\infty}\leq M$. Let $u$ be the unique weak solution  of \eqref{equi} in $L^{\infty}_{loc}\left([0, T_{max}(u_0)); W^{1,\infty}(\Omega)\right)$. 
Then
\begin{equation}\label{problo}
 |\nabla u|\leq C_1 \delta^{-1/(q-p+1)} (x)+C_2.\qquad\text{in} \quad\Omega\times \left(0, T_{max}(u_0)\right).
\end{equation}
 where $C_1=C_1(q,p,N)>0$ and $C_2=C_2(q,p,\Omega, M, \left\|g\right\|_{C^2})>0$.
\end{theor}
This estimate in particular implies that $|\nabla u|$ remains bounded away from the boundary. Therefore, when $T_{max} (u_0)<\infty$,  the blow-up may only take place on the boundary and \eqref{problo} provides information on the blow-up profile near $\partial\Omega$. Estimate \eqref{problo} is sharp in one space dimension, see \cite{amal}. Similar results are already available for $p = 2$ and have been established in \cite{Zhan},\cite{Arri}. For $p>2$, only global-in-space gradient estimates were available up to now (ie for $\Omega=\mathbb{R}^N$, see \cite{bartier}).  The proof of  estimate \eqref{problo} is  based  on similar arguments as for the case $p = 2$, namely Bernstein type arguments, but they are much more technical.
Moreover, the proof of (\ref{problo}) also relies on  a regularizing effect for solutions to (\ref{equi}) which seems to be new and which is stated below.
\begin{theor}\label{rek}
 Assume that $q >  p-1>1$ and let  $u$ be the unique weak solution of problem (\ref{equi}) in $ L^{\infty}_{loc}\left([0,T_{max}(u_0));  W^{1,\infty}(\Omega)\right)$. Then
\begin{equation}
u_t\leq \dfrac{1}{p-2} \dfrac{\left\|u_0\right\|_{\infty}}{t}  \ \ \ \ \mathrm{in} \  \ \mathcal{D}'(\Omega) \ \ \mathrm{a,e} \ \, t>0.
\end{equation}
\end{theor}
Let us note that due to the positivity of the source term, this inequality implies the semi-concavity estimate 
\begin{equation}
\Delta_p (u)= \text{div}\left( |\nabla u|^{p-2} \nabla u\right) \leq \dfrac{C}{t},
\end{equation}
which was obtained in the case $\Omega= \mathbb{R}^N$ by a different method in \cite{este}.

Finally we give the following blow-up result, which is a variant of a global nonexistence result in \cite{oups}, reinterpreted in terms of GBU in the light of Theorem \ref{exis}. Let $\varphi_1$ be the first eigenfunction of $-\Delta$ with homogeneous Dirichlet boundary conditions
\begin{theor}\label{blo}
Assume that \textbf{$q>p>2$} and let $u$ be the unique weak solution of (\ref{equi}) in $L^{\infty}_{loc}\left([0,T_{max}(u_0));  W^{1,\infty}(\Omega)\right)$. Let $\alpha \geq 1$ satisfy $\dfrac{p-1}{q-p+1}<\alpha< q-1$, then there exists a constant $C=C(q,p,\alpha,\Omega,\left\|g\right\|_{\infty})>0$ such that if
$\int_{\Omega} u_0 \, \varphi_1^{\alpha} dx\geq C$, then $T_{max}(u_0)<\infty$, i.e.  gradient blow-up occurs.
\end{theor}
For results concerning other aspects of equation \eqref{equi} and the corresponding Cauchy problem, see e.g. \cite{chen,shi,chen2, zhao2,bartier} and the references therein. Asymptotic behavior of global solution is investigated in \cite{stin,barle,oups,lauren1, lauren3,laurenew, andreu} and references therein.

The rest of the paper is organized as follows: In Section 2 we prove the well-posedness  of (\ref{equi}) in $W^{1,\infty}(\Omega)$, as well as the regularizing effect.
Section 3 is devoted to the proof of Theorem \ref{profil}.
Finally in section 4 we prove the sufficient blow-up criterion of Theorem \ref{blo}.
\section{Proof of Theorem \ref{exis} and Theorem \ref{rek}}
\subsection{Local existence}

Consider the following approximate problems for \eqref{equi}:
\begin{equation}\label{apro}
\left\{
\begin{array}{lll}
 \partial_t u_n- \mathrm{div} \left(\left( |\nabla u_n|^2+ \dfrac{1}{n}\right)^{(p-2)/2} \nabla u_n \right)= \left( |\nabla u_n|^2+ \dfrac{1}{n} \right)^{q/2}-\dfrac{1}{n^{q/2}},  & x\in \Omega, t>0,\\
u_n(x,t)=g(x),    & x \in\partial\Omega, t>0,\\
u_n(x,0)=u_0(x),      & x\in\Omega.
\end{array}
\right.
\end{equation}
For each fixed $n\in\mathbb{N}$, problem \eqref{apro} is  no longer degenerate and the regularity theory of quasilinear parabolic equations \cite{lady} provides  local-in-time solutions  $u_n$, which are smooth for $t>0$ and continuous up to $t=0$.

To find the limit function $u(x,t)$ of the sequence $\left\{ u_n(x,t)\right\}$, we divide our proof into  5 steps. Recall that there exists $\eta_0> 0$ small such that, for  any $x\in\overline{\Omega}$ with $\delta(x)\leq \eta_0$, the point $\widetilde{x}:=\text{proj}_{\partial\Omega} (x)$ (the projection of $x$ onto the boundary) is well defined and unique.\\

\textbf{STEP 1.} There exist a small time $T_0>0$, $\eta\in (0, \eta_0)$ and $ M_2>0$, all independent of $n$ and depending on  $u_0$  through $M$ only, such that 
\begin{equation}\label{max}
\left\|u_n\right\|_{L^{\infty}(Q_{T_0})}\leq M_1:= \text{max}\left( \left\|u_0\right\|_{\infty},\left\|g\right\|_{\infty}\right),
\end{equation}
and
\begin{equation}\label{sert}
\underset{\substack{x\in\Omega\\
\delta(x)\leq\eta}}{\text{sup}} \dfrac{| u_n(x,t)-u_n(\widetilde{x},t)|}{\delta(x)}\leq M_2, \qquad 0<t\leq T_0.
\end{equation}

The barrier function will have the form
Estimate \eqref{max} is a direct consequence of the maximum principle since $M_1$ is a super solution for any $n$.\\

In order to prove estimate \eqref{sert}, we are going to construct a local barrier function under the exterior sphere condition satisfied by the domain $\Omega$, i.e. for any $x$ near $\partial\Omega$, a supersolution in a neighborhood of $x$.
  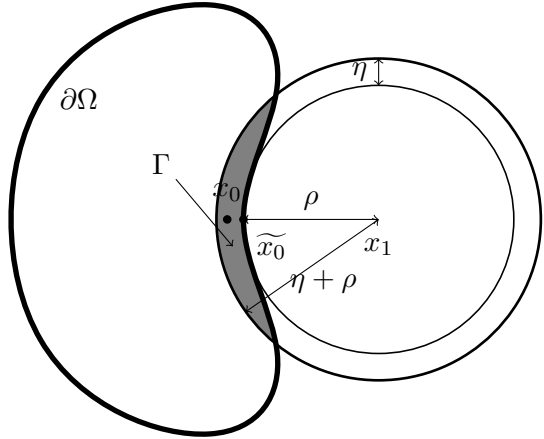
\begin{figure}[!h!]
	\centering
		\input{mondessin.pgf}
		\caption{Local barrier function}
	\label{mondessin}
\end{figure}

Let $\rho>0$ be such that for all $x\in\partial\Omega$, $\overline{B_{\rho}(x+\rho\nu_x)}\cap\overline{\Omega}=\left\{x\right \}$, where $\nu_{x}$ is the unit outward normal vector on $\partial\Omega$ at $x$. Fix an arbitrary  $x_0\in\Omega$ such that $\delta(x_0)\leq \eta$ where $\eta\in (0, \eta_0)$ will be chosen later. Define $x_1=\widetilde{x_0}+\rho\nu_{\widetilde{x_0}}$.  Without loss of generality we may assume that $x_1=0$ and we write $r=|x|$. Let us denote, for $s\geq 0$, 
 \begin{equation}\label{jil}
 a(s)=\left(s+\dfrac{1}{n}\right)^{(p-2)/2}, \quad\text{and}\qquad \kappa=\dfrac{2 a'(s) s}{a(s)}\in [0, p-2].
 \end{equation}
 
We recall that for a function $\phi(x)=\phi(|x|)$, we have:
 \begin{equation}\label{rapp}
 \begin{array}{ll}
& \nabla \phi(x)=\phi'(r) \dfrac{x}{r},\\
&D^2 \phi(x)=\phi''(r) \dfrac{x\otimes x}{r^2}+ \dfrac{\phi'(r)\text{Id}}{r} -\phi'(r) \dfrac{x    \otimes x}{r^3},\\
&\Delta \phi(x)=\phi''(r)+\dfrac{(N-1)\phi'(r)}{r},
\end{array}
\end{equation}
 where $\text{Id}$ is the unit matrix and $(x\otimes x)_{ij}=x_i x_j$.

$$\bar{v}(x,t)= \phi (r-\rho) + g(x),$$
where $\phi$ is a smooth function of one variable which is increasing and \textbf{concave}.
First let us write 
\begin{align}\label{expres} \mathrm{div}\left(\left(|\nabla\bar{v}|^2+\dfrac{1}{n}\right)^{(p-2)/2}\nabla\bar{v}\right)&=a(|\nabla \bar{v}|^2)\Delta \bar{v}+2a'(|\nabla\bar{v}|^2) (\nabla\bar{v})^t D^2\bar{v} \nabla \bar{v},\nonumber\nonumber\\
&=a(|\nabla\bar{v}|^2)\left( \Delta \bar{v}+\kappa \,\frac{(\nabla\bar{v})^t D^2\bar{v} \nabla \bar{v}}{|\nabla \bar{v}|^2}\right).
\end{align}
Using \eqref{rapp}, we have 

\begin{eqnarray*}
\begin{array}{lll}
&\left [ \Delta \bar{v}+\kappa\,\dfrac{(\nabla\bar{v})^t D^2\bar{v} \nabla \bar{v}}{|\nabla \bar{v}|^2}\right]\\
&= \phi''(r-\rho)+\dfrac{(N-1)\phi'(r-\rho)}{r}+\Delta g \\
 &+\kappa \,\dfrac{\phi''(r-\rho) (\nabla\bar{v}\cdot x)^2}{r^2|\nabla\bar{v}|^2}+ \kappa \,\dfrac{\phi'(r-\rho)}{r}-\kappa\, \dfrac{\phi'(r-\rho) (\nabla\bar{v}\cdot x)^2}{r^3|\nabla\bar{v}|^2} +\kappa\, \dfrac{(\nabla\bar{v})^t D^2  g  \nabla\bar{v}}{|\nabla\bar{v}|^2}.
\end{array}
\end{eqnarray*}
Since $\phi'(r-\rho)\geq 0$, $r\geq \rho$, $\kappa\geq 0$ and  $0\geq\phi''(r-\rho)$, we have

\begin{equation}\label{cop}
-\left[ \Delta \bar{v} +\kappa\,\frac{(\nabla\bar{v})^t D^2\bar{v} \nabla \bar{v}}{|\nabla \bar{v}|^2}\right]\geq -\phi''(r-\rho)-\dfrac{(N-1+\kappa)}{\rho} \phi'(r-\rho) -\left\|\Delta g\right\|_{\infty}-\kappa\,\left\|D^2 g\right\|_{\infty}.
\end{equation}

On the other hand $|\nabla\bar{v}|=\left|\phi'(r-\rho) \dfrac{x}{r}+ \nabla  g\right|\leq \phi'(r-\rho)+|\nabla g |\leq 2 \phi'(r-\rho)$ provided that 
\begin{equation}\label{jal}
\phi'(r-\rho)\geq  \left\|\nabla g \right\|.
\end{equation}
In this case we have 
\begin{equation}\label{cip}
\left(|\nabla\bar{v}|^2+\frac{1}{n}\right)^{(q-p+2)/2}\leq \big[4(\phi'(r-\rho))^2+1\big]^{(q-p+2)/2}.
\end{equation}
We take
\begin{equation*}
\phi(s)=s(s+\delta)^{-\beta}, \qquad s\geq 0,
\end{equation*}
 where $\beta=\beta(q,p)\in(0,1)$ is to be chosen later. We denote $\Gamma := B(x_1, \rho+\eta)\cap\Omega$ (see figure \ref{mondessin}).  Our aim is to show that $\bar{v}$ is a super-solution in $\Gamma\times(0,T_0)$ where $T_0, \delta>0$ and $\eta \in (0, \eta_0)$ small enough. In the rest of the proof, the constants $T_0, \eta, \delta$ and $C$ will be independent of $x_0$, $n$ and will depend on the initial data $u_0$ through $M$ only (and they will depend on the other data $p, q, N, \Omega$ and $\left\|g\right\|_{C^2}$ without other mention).
We calculate 
\begin{align*}
& \phi'(s)=\left[ (1-\beta)s+\delta  \right] \left( s +\delta\right)^{-\beta-1},\\
&\phi''(s)=-\beta \left[(1-\beta)s+2 \delta\right] \left( s +\delta\right)^{-\beta-2}.
\end{align*}We are looking for condition on $\beta$ and $\delta$ such that
\begin{equation}\label{ret}
-\mathrm{div}\left(\left(|\nabla \bar{v}|^2+\dfrac{1}{n}\right)\nabla \bar{v}\right)\geq \left( |\nabla\bar{v}|^2+\frac{1}{n} \right)^{q/2}-\left(\frac{1}{n}\right)^{q/2}.
\end{equation}
Due to \eqref{expres}, it suffices to have
\begin{equation}
-\left[ \Delta \bar{v} +\kappa\,\frac{(\nabla\bar{v})^t D^2\bar{v} \nabla \bar{v}}{|\nabla \bar{v}|^2}\right]\geq\left( |\nabla \bar{v}|^2+\dfrac{1}{n}\right)^{\frac{q-p+2}{2}},
\end{equation}
which, by \eqref{cop}-\eqref{jil}-\eqref{cip} reduces to 
\begin{equation}
-\phi''(r-\rho)+\dfrac{(3-N-p)}{\rho} \phi'(r-\rho) \geq \big[4(\phi'(r-\rho))^2+1\big]^{(q-p+2)/2}+ (p-2+\sqrt{N})\left\|D^2 g\right\|_{\infty}.
\end{equation}
Therefore (\ref{ret})  holds if
\begin{align*}
 (r-\rho+\delta)^{-\beta-2}\left[2\beta\delta+ (3-N-p) \dfrac{(\eta+\delta)^2}{\rho}\right]
&\geq \left[ 4(r-\rho+\delta)^{-2\beta}+1\right]^{(q-p+2)/2} \\
&+(p-2+ \sqrt{N})\left\|D^2  g \right\|_{\infty}.
\end{align*}
Assume that  $\eta$ and  $\delta$ are  such that 
\begin{eqnarray}
\left\{
\begin{array}{lll}
4(r-\rho+\delta)^{-2\beta}\geq 4(\eta+\delta)^{-2\beta}\geq1,\\
2\beta\delta+\dfrac{(3-N-p)}{\rho} (\eta+\delta)^2\geq \beta\delta,\\

\end{array}
\right.
\end{eqnarray}
then to get (\ref{ret}) it is sufficient to have 
\begin{equation}\label{ing}
\beta\delta(r-\rho +\delta)^{-\beta-2}\geq(r-\rho+\delta)^{-\beta(q-p+2)}4^{(q-p+3)},
\end{equation}
and
\begin{equation}
\beta\delta(r-\rho +\delta)^{-\beta-2}\geq 4 (p-2 + \sqrt{N})\left\|D^2 g \right\|_{\infty}.\label{inegs}
\end{equation}
Inequality \eqref{ing} holds if we choose $\eta=\delta, \, \beta=\frac{1}{2(q-p+2)},$ and $\delta$ satisfying
$$4^{p-q-4} \beta\geq  \delta^{\frac{q-p+3}{2(q-p+2)}}.$$
Inequalities \eqref{ing}-\eqref{inegs} and \eqref{jal} hold if we choose $\delta$ small enough. We have thus shown that if $\eta=\delta$ is small, then $\bar{v}$ is a supersolution on $\Gamma\times(0, T_0)$ for any $T_0>0$.

Now we need to have a control on the parabolic boundary of  $\Gamma\times(0,T_0)$ for $T_0> 0$ small. For this purpose,  we introduce another comparison function
$$\bar{u}(x,t)=(C^2K^2+1)^{q/2}t+ C(1-e^{-K(r-\rho)})+\left\| g \right\|_{\infty}.$$ It is easy to see that  if we choose $K$ sufficiently large $\left(K>\dfrac{N+p-3}{\rho}\right)$, then
$$-\mathrm{div}\left(\left(|\nabla\bar{u}|^2+\frac{1}{n}\right)^{(p-2)/2} \nabla \bar{u}\right)\geq 0.$$
Thus $$\partial_t \bar{u}-\mathrm{div}\left(\left(|\nabla\bar{u}|^2+\frac{1}{n}\right)^{(p-2)/2} \nabla \bar{u}\right)\geq \left( |\nabla \bar{u}|^2+ \dfrac{1}{n} \right)^{q/2}-\left(\dfrac{1}{n}\right)^{q/2}.$$
 Next we may choose $C>0$ large enough (depending only on $M$) such that $C(1-e^{-K(r-\rho)})+\left\| g \right\|_{\infty}\geq u_0(x)$ in $\Omega$. Since $\bar{u}\geq g$ on $\partial\Omega\subset\left\{ x\in\mathbb{R}^N, \quad |x|\geq \rho\right\}$, by the maximum principle we get that for any $n$, $ u_n\leq \bar{u}$ in $\Omega_T$.
Thus  \begin{align*}
 u_n(x,t)&\leq  (C^2 K^2+1)^{q/2} t+ C(1-e^{-K\eta})+\left\| g \right\|_{\infty}\\
 & \leq 2^{-\beta} \eta^{1-\beta}+ g (x)=\bar{v}(x,t)
 \end{align*}
  on $\left\{x\in\Omega, |x|=\rho+\eta\right\}\times[0, T_0]$, provided   $T_0$ and $\eta=\delta$ are small enough (depending only on $M,p,q, \Omega, \left\|g\right\|_{C^2}$).

On the other hand $u=g\leq \bar{v}$ on $\partial\Omega\times[0,T_0]$. We conclude that $\bar{v}$ is a super solution on $\Gamma\times (0, T_0)$. Similarly $\underline{v}:= g-\phi(r-\rho)$ is a sub-solution. Applying the maximum principle we get $\underline{v}\leq u_n\leq\bar{v}$ on $\Gamma\times[0,T_0]$, and hence in particular
$$\dfrac{|u_n(x_0,t)-u_n(\tilde{x_0},t)|}{|x_0-\tilde{x_0}|}\leq \underset{0\leq s\leq \delta}{\text{sup}} |\phi'(s)|+ \left\|\nabla g\right\|_{\infty}\leq \delta^{-\beta} +\left\|\nabla g\right\|_{\infty}=: M_2, \quad 0<t\leq T_0,$$
which yields \eqref{sert}.\\

\textbf{STEP\,2}. There holds 

\begin{equation}\label{grt}
\left\|\nabla u_n\right\|_{L^{\infty}(Q_{T_0})}\leq M_3:= M_2+\left\|\nabla g\right\|_{\infty}.
\end{equation}
We use a similar argument as in \cite[Theorem 5]{KK}. Let $h\in\mathbb{R}^N$ satisfy $|h|\leq \eta$. Due to the translation invariance of \eqref{apro}, if $u_n$ is a classical solution of (\ref{apro}) in $\Omega$, then the function $u^h_n:= u_n(x-h,t)$ is a classical solution of (\ref{apro}) in  $\Omega_h\times (0,T_0)$ where $\Omega_h:=\left\{x\in \mathbb{R}^N \, |\ \  x-h\in\Omega\right\}$.
 Let $t\in [0, T_0]$ and $x\in\partial(\Omega\cap\Omega_h)$. We may assume for instance $x\in\partial\Omega$, the case $x+h\in \partial\Omega$ being similar. Then using $|\widetilde{y}-\widetilde{z}|\leq |y-z|$ and \eqref{sert}, we get
\begin{eqnarray*}
|u_n(x,t)- u_n(x+h,t)|&=|u_n(\widetilde{x},t)- u_n(\widetilde{x+h},t)+u_n(\widetilde{x+h},t)- u_n(x+h,t)|\\
&\leq \left\|\nabla g\right\|_{\infty} |\widetilde{x}-\widetilde{x+h}|+M_2 \delta(x+h)\\
&\leq (\left\|\nabla g\right\|_{\infty}+M_2)|h|= M_3|h|.
\end{eqnarray*}
In particular $u_n(x,t)\leq u_n^h(x,t)+M_3|h|$ on $\partial (\Omega\cap\Omega_h)\times [0, T_0]$.

Applying the maximum principle, we have $u_n(x,t)\leq u_n^h(x,t)+M_3|h|$ on  $(\Omega\cap\Omega_h)\times [0, T_0]$. By the same argument $u_n^h(x,t)-M_3|h|\leq u_n(x,t)$ on  $(\Omega\cap\Omega_h)\times [0, T_0]$, hence $|u_n(x,t)-u_n^h(x,t)|\leq M_3 |h|$. Since $|h|\leq\eta$ is arbitrary, the conclusion follows.\\

\textbf{STEP\ 3}. Let $\epsilon>0$ and set  $Q_{T_0,\epsilon}=\left\{x\in\Omega, \delta(x)>\epsilon\right\}\times(\epsilon, T_0-\epsilon)$. There exists a constant $M_4>0$ independent of $n$, such that 
\begin{equation}\label{bened}
\left|\nabla u_n(x_1, t_1)-\nabla u_n(x_2,t_2)\right|\leq M_4 \left(|x_1-x_2|^{\alpha}+|t_1-t_2|^{\frac{\alpha}{2}}\right) 
\end{equation}
for any pair of points $(x_i,t_i)\in Q_{ T_0 ,\epsilon}$, where $M_4$ and $\alpha$ are positive constants depending only on $ T_0, M_3$ and  $\epsilon$. Indeed  we know from a result of DiBenedetto and Friedman \cite{dib} that if $f\in L^r(\Omega_T)$ for some $r> \dfrac{pN}{p-1}$ then weak solutions of degenerate parabolic equation of the form 
\begin{equation}
 \partial_t v -\mathrm{div} \left( |\nabla v|^{p-2}\nabla v\right)=f(x,t)
\end{equation}
are of class $C^{1,\alpha}_{loc}(Q_T)$ with $\mathrm{H\ddot{o}lder}$ norm depending only on $\left\|f\right\|_{L^r}, \left\|\nabla u\right\|_{L^p}$ and  $\left\|u\right\|_{L^{\infty}_t,L^2_x}$.\\

\textbf{STEP\,4}. There exists a constant $M_5>0$ independent of $n$, such that 
\begin{equation}\label{time}
\left\|\partial_t u_n\right\|_{L^2(Q_{T_0})}\leq M_5.
\end{equation}
To see this , multiplying \eqref{apro} by $\partial_t u_n$ and integrating over $Q_{T_0}$, we have
\begin{align*}
\int_0^{T_0}\int_{\Omega}(\partial_t u_n)^2 dx dt =&- \int_0^{T_0}\int_{\Omega} \left( |\nabla u_n|^2+\dfrac{1}{n}\right)^{(p-2)/2} \nabla u_n\cdot\nabla (\partial_t u_n) dx dt\\
&+\int_0^{T_0}\int_{\Omega}\left( |\nabla u_n|^2+\dfrac{1}{n}\right)^{q/2} \, \partial_t u_n\, dx dt.
\end{align*}
By $\mathrm{H\ddot{o}lder}$'s inequality and
\begin{align*}
&\int_0^{T_0}\int_{\Omega} \left( |\nabla u_n|^2+\dfrac{1}{n}\right)^{(p-2)/2} \nabla u_n\cdot\nabla (\partial_t u_n) dx dt\\
&=\frac{1}{p}\int_{\Omega} \left(|\nabla u_n(x,T_0)|^2+\dfrac{1}{n}\right)^{p/2}-\frac{1}{p}\int_{\Omega}\left(|\nabla u_n(x,0)|^2+\dfrac{1}{n}\right)^{p/2},
\end{align*}
we get
\begin{align*}
\int_0^{T_0}\int_{\Omega}(\partial_t u_n)^2 dx dt &\leq\frac{2}{p}\int_{\Omega}\left( |\nabla u_n|(x,0)|^2+\dfrac{1}{n}\right)^{p/2}dx 
+ 2 \int_0^{T_0}\int_{\Omega}\left( |\nabla u_n|^2+\dfrac{1}{n}\right)^{q} dx dt\\
&\leq M'.
\end{align*}
 for some $M'=M'\big( |\Omega|, M_3, T_0, p, q\big)>0$.
\newline

\textbf{STEP\,5}.
We recall that by the Rellich-Kondrachev theorem we have
\begin{equation}\label{compact}
W^{1, \infty}(\Omega)\stackrel{c}\hookrightarrow C(\overline{\Omega})\hookrightarrow\ L^2(\Omega).
\end{equation}
Using (\ref{max})-(\ref{grt})-(\ref{time})-\eqref{compact} and the compactness theorem in  [\cite{simon} Corollary 4], we have that $\left\{u_n\right\}$ is relatively compact in $C\left( [0,T_0]; C(\overline{\Omega})\right)=C\left(\overline{\Omega}\times[0, T_0]\right)$. 
 By virtue of (\ref{grt})-(\ref{bened})-(\ref{time}), the Ascoli-$\mathrm{Arzel\grave{a}}$ theorem and the relative compactness of $\left\{u_n\right\}$ in $C\left(\overline{\Omega}\times[0, T_0]\right)$, we can find a subsequence, still denoted by $\left\{ u_n\right\}$ for convenience, such that, for each $\epsilon>0$, 
\begin{eqnarray}\label{fin}
\left.
\begin{array}{lll}
u_n\to u \,\,\,\, &\quad\text{in}\quad  C\left(\overline{\Omega}\times [0,T_0]\right),\\ 
\nabla u_n\to \nabla u \quad &\quad\text{in}\quad  C(Q_{T_0,\epsilon} ),\\
\partial_t u_n\to \partial_t u \quad &\text{weakly in}\quad  L^2(Q_{T_0}).
\end{array}
\right\}
\end{eqnarray}
 We multiply \eqref{apro} by a test function and integrate. Then by the Lebesgue's dominated convergence theorem and \eqref{fin} we can pass to the limit and check that $u$ is a weak solution of \eqref{equi}.

\subsection{The blow-up alternative}

Let us temporarily assume the uniqueness result  which will be proved in the next section. The construction of the weak solution as a limit of classical solutions  implies the blow-up alternative.

Indeed suppose that the maximal existence time $T_{max}(u_0)<\infty$ and that there exist $\mathcal{M}>0$ and $t_k\to T_{max}(u_0)$ such that for all $k$ 
\begin{equation}\label{bv}
 \left\|\nabla u(t_k)\right\|_{L^{\infty}(\Omega)}\leq \mathcal{M}.
\end{equation}
Then we can find $\tau=\tau(\mathcal{M})>0$ independent of $k$, such that  the problem

\begin{equation}
 \left\{
\begin{array}{lll}
 \partial_t u- \mathrm{div}( |\nabla u|^{p-2} \nabla u)=|\nabla u|^q,  &  x\in\Omega, t>0,\\
u(x,t)=g(x),                                                  & x\in\partial\Omega, t>0,\\
u(x,0)=u(x, t_k),                                             & x\in\Omega,
\end{array}
\right.
\end{equation}
admits a unique weak solution $v_k$ on $[0, \tau)$. Setting  $\tilde{u}(t)=\left\{\begin{array}{ll}
u(t) &\text{for} \, t\in[0, t_k)\\
v_k(t-t_k) &\text{for} \, t\in[t_k, t_k+\tau)
\end{array}
\right.$,  it is easy to see that we get  a weak solution defined on $[0, t_k+\tau)$.

Since for $k$ large enough $t_k+\tau > T_{max}(u_0)$, this contradicts the definition of $T_{max}(u_0)$. Hence $T_{max}(u_0)<\infty\Rightarrow\underset{t\to T_{max}(u_0)}{\lim} \left\|\nabla u\right\|_{L^{\infty}(\Omega)}=\infty$.

\subsection{Uniqueness}

In this section we prove the uniqueness of the weak solution. This result will be a consequence of the following comparison principle which, in turns, also guarantees \eqref{pmax}.
\begin{propo}\label{comprin}
Let $u,v$  be respectively, sub-, super-solutions of \eqref{equi}. Assume that $u, v\in  L^\infty\left((0,T);W^{1,\infty}(\Omega)\right)$.
Then $u\leq v$ on $ \Omega \times (0,T)$.
\end{propo}

The proof of  Proposition  \ref{comprin} is mostly based on the following algebraic
lemma from which we can show that the source term can be
counterbalanced by the diffusion effect (c.f \cite{cool}).	
\begin{lem}[Monotonicity Property]\label{lemmefin}
Let $\sigma>1$. For all  $a$ and  $b$ $\in\mathbb{R}^N$:

$$\left\langle  |a|^{\sigma-2} a-|b|^{\sigma-2} b, a-b\right\rangle \geq 
\dfrac{4}{\sigma^2} \Big| |a|^{(\sigma-2)/2} a- |b|^{(\sigma-2)/2} b \Big|^2.
 $$
\end{lem}
\textbf{Proof of Proposition \ref{comprin}}. We set  $w=(u-v)^+$. By definition we have $w=0$ on $\partial\Omega$.
By Remark \ref{rem1}, for any $\tau \in (0, T)$, using $\psi=w$ as  test-function, we have
\begin{align*}
&\int_0^{\tau} \int_{\Omega} w w_t \,dx dt\leq\\
&\underbrace{\int_0^{\tau}\int_{\left\{w(\cdot, t)>0\right\}} \left[ |\nabla u|^q-|\nabla v|^q \right] w \,dx dt}_{\mathcal{B}} -\underbrace{\int_0^{\tau} \int_{\left\{w(\cdot,t)>0\right\}} \left[ |\nabla u|^{p-2} \nabla u- |\nabla v|^{p-2} \nabla v\right]\cdot \nabla w \,dx dt}_{\mathcal{H}}.
\end{align*}
We set $a=\nabla u$ and $b=\nabla v$. We get by lemma \ref{lemmefin}
\begin{eqnarray}\label{bona}
\mathcal{H}\geq c(p)\int_0^{\tau} \int_{\left\{w(\cdot, t)>0\right\}} \left| |\nabla u|^{(p-2)/2} \nabla u- |\nabla v|^{(p-2)/2} \nabla v \right|^2 dx dt.
\end{eqnarray}
Let's consider the term  $\mathcal{B}$. We put  $h(s) =s^{\frac{2q}{p}}$ for $s \geq 0$. Given that $q\geq p-1\geq\dfrac{p}{2}$, we have  $h'(s) = \frac{2q}{p} s^{\frac{2q-p}{p}}$. The mean value theorem yields
$$\Big| |\nabla u|^q-|\nabla v|^q  \Big|^2\leq C h'(\theta) ^2 \left| |\nabla u|^{(p-2)/2} \nabla u- |\nabla v|^{(p-2)/2}\nabla v \right|^2,$$
for some $0\leq \theta \leq \mathrm{max}(|\nabla u|^{\frac{p}{2}}, |\nabla v|^{\frac{p}{2}})$.\\
Since we assumed $u,v\in L^\infty\left((0,T);W^{1,\infty}(\Omega)\right)$, it follows that

$$ \Big| \, |\nabla u|^{q }-|\nabla v|^q \, \, \Big|^2\leq
C\left| |\nabla u|^{(p-2)/2}\nabla u-|\nabla v|^{(p-2)/2}\nabla v \right|^2.$$
On the other hand, the Young inequality implies
$$\mathcal{B}\leq \epsilon \int_0^{\tau}\int_{\left\{w(\cdot, t)>0\right\}} \left| \,|\nabla u|^{q}-|\nabla v|^{q}\, \right|^2 dx dt+ C(\epsilon)\int_0^{\tau} \int_{\left\{w(\cdot, t)>0\right\}} w^2 dx dt.$$
Combining these two inequalities, we arrive at
\begin{equation}\label{bon}
\mathcal{B}\leq  C\epsilon \int_0^{\tau}\int_{\left\{w(\cdot,t)>0\right\}} \left| |\nabla u|^{(p-2)/2} \nabla u-|\nabla v|^{(p-2)/2} \nabla v \right|^2 dx dt+ C(\epsilon) \int_0^{\tau}\int_{\left\{w(\cdot, t)>0\right\}} w^2 dx dt.
\end{equation}
Choosing $\epsilon$ small enough, we get
\begin{equation}
\int_{\Omega} w^2(\tau) \,dx \leq\int_{\Omega} w^2(0) \,dx + C(\epsilon) \int_0^{\tau}\int_{\Omega} w^2 \,dx dt, \quad 0<\tau< T.
\end{equation}
The Gronwall lemma implies that for any $t\in (0,T)$
$$\int_{\Omega} w^2(x,t)\ dx\leq  e^{C t } \int_{\Omega} w(x,0)^2 \ dx .$$
We conclude that $w\equiv0$ almost everywhere.

\begin{rem}
\begin{enumerate}[(a)]
 \item 
The inequality in lemma \ref{lemmefin}  for $\sigma\in(1,2)$ can be deduced from the inequality for  $\sigma\geq 2$ in \cite{cool} as follows:\\

We set  $a=|\nabla u|^{\sigma-2}\nabla u$ and $b=|\nabla v|^{\sigma-2} \nabla v$.
\begin{eqnarray}
\begin{array}{ll}
\left\langle |\nabla u|^{\sigma-2}\nabla u- |\nabla v|^{\sigma-2}\nabla v, \,\nabla u-\nabla v\right\rangle&=\left\langle a-b, a \,|a|^{\frac{2-\sigma}{\sigma-1}} -\,b \,|b|^{\frac{2-\sigma}{\sigma-1}}\right\rangle\\
&=\left\langle a-b,\, a\,|a|^{m-2}- b\, |b|^{m-2}\right\rangle.
\end{array}
\end{eqnarray}
where $m= \frac{\sigma}{\sigma-1} > 2$.

\item
The question of uniqueness was partially open in \cite{stin}. The  preceding result can be applied to show uniqueness in the case $p-1\geq q\geq \dfrac{p}{2}$  with $p\geq 2$.

\item
In \cite{anton} we have a weaker inequality for $p\in (1, 2)$ but it is sufficient to prove uniqueness  for the case $q>1$:
$$\left\langle  |a|^{p-2} a-|b|^{p-2} b, a-b\right\rangle \geq (p-1) |a-b|^2\left(|a|^{p}+|b|^{p}\right)^{\frac{p-2}{p}}.$$
\end{enumerate}

\end{rem}

\subsection{Regularizing effect}

We use a technique developed by Zhao for the
the $p$-Laplace equation without source term \cite{zhao}. The idea is to apply a
Stampacchia maximum principle argument to the equation satisfied by $\lambda^{\gamma} u(x,\lambda t)- u(x,t)$
 and then let $\lambda\to 1^+$.
Let $u$ be a weak solution of \eqref{equi} in $L^{\infty}_{loc}\left([0,T); W^{1,\infty}(\Omega)\right)$. Set 
$$u_{\lambda}(x,t)=\lambda^{\gamma} u(x,\lambda t), \qquad\lambda>1, \gamma=\frac{1}{p-2}.$$
Then $u_{\lambda}$ is a weak solution of
\begin{equation*}
 \left\{
\begin{array}{lll}
 \partial_t u_{\lambda}- \mathrm{div}( |\nabla u_{\lambda}|^{p-2} \nabla u_{\lambda})=\lambda^{-(q-p+1)\gamma}|\nabla u_{\lambda}|^q,  \,\,\,\,\,&  x\in\Omega, t\in\left(0,\frac{T}{\lambda}\right),\\
u_{\lambda}(x,t)=\lambda^{\gamma} g(x),                                                  & x\in\partial\Omega,  t\in\left(0, \frac{T}{\lambda}\right),\\
u_{\lambda}(x,0)=\lambda^{\gamma} u_0(x),                            & x\in\Omega.

\end{array}
\right.
\end{equation*}
Let $w=u_{\lambda}-u$. Then using Remark \ref{rem1}, we have for any $\phi\in C^0(\overline{Q_T})\cap L^p((0,T); W^{1,p}(\Omega))$, $\phi=0$ on $\partial\Omega\times(0, T)$ and  any $0<\tau<\dfrac{T}{\lambda}$,
\begin{align*}
&\int_0^{\tau}\int_{\Omega} \phi  w_t \,dx dt=\\
&\int_0^{\tau}\int_{\Omega} \left[ \lambda^{-(q-p+1)\gamma}|\nabla u_{\lambda}|^q-|\nabla u|^q \right] \phi \,dx dt -\int_0^{\tau}\int_{\Omega} \Big[ |\nabla u_{\lambda}|^{p-2} \nabla u_{\lambda} - |\nabla u|^{p-2} \nabla u\Big]\cdot \nabla\phi  \,dx dt.
\end{align*}
Taking
$$\phi= (w-k)_{+}, \qquad k=(\lambda^{ \gamma}-1)\left\| u_0 \right\|_{L^{\infty}(\Omega)};$$
we have $\phi=0$ on $\partial\Omega$ due to \eqref{initial}. Since  $\lambda>1$ and $\phi\geq0$, we get

\begin{align*}
\int_0^{\tau}\int_{\Omega} (w-k)_+\partial_t (w-k)_+& \,dx  dt\leq \int_0^{\tau}\int_{\left\{w(\cdot, t) > k\right\}} \left[ |\nabla u_{\lambda}|^q-|\nabla u|^q \right] (w-k)_+ \,dx dt\\
&-\int_0^{\tau}\int_{\left\{w(\cdot, t) > k\right\}} \left[ |\nabla u_{\lambda}|^{p-2} \nabla u_{\lambda} - |\nabla u|^{p-2} \nabla u\right]\cdot (\nabla u_{\lambda}-\nabla u) \,dx dt.
\end{align*}
Using the same trick as for (\ref{bona})-(\ref{bon}), we get 
$$\int_{\Omega} (w-k)_+^2(x,t) dx\leq \left(\int_{\Omega}(w(x,0)-k)_+^2 dx\right) e^{Ct}.$$
Given that  $(\lambda^{\gamma}-1)u_0(x)\leq (\lambda^{\gamma}-1)\left\|u_0\right\|_{L^{\infty}}$, we get $(w-k)_+ \equiv 0$ a.e on $\Omega\times \left(0,\frac{T}{\lambda}\right)$. Thus 
\begin{equation}\label{serk}
 \lambda^{\gamma} u(x,\lambda t)-u(x,t)\leq (\lambda^{\gamma}-1)\left\|u_0\right\|_{L^{\infty}}.
\end{equation}
Dividing (\ref{serk}) by $(\lambda-1)$ and letting $\lambda\to1^+$, we get 
$$\gamma u(x,t)+ t u_t(x,t)\leq\gamma \left\|u_0\right\|_{L^{\infty}}.$$
We conclude using the positivity of $u$.

\begin{rem}
 The homogeneity of the operator and the boundedness of $u$ are essential.
\end{rem}

\section{Gradient estimate: proof of Theorem \ref{profil}}

The proof of \eqref{problo} relies on a modification of the Bernstein technique and  the use of a suitable cut-off function. It requires the study of the partial differential equation satisfied by $|\nabla u|^2$. We  follow the  ideas used in \cite{Zhan} and \cite{bartier}.
Let $x_0\in\Omega$ be fixed, $0< t_0<T<T_{max}(u_0)$, $R>0$ such that $B(x_0, R)\subset\Omega$ and write $Q^{t_0}_{T, R}= B(x_0, R)\times (t_0,T)$

Let $\alpha \in  (0,1)$ and set $ R' =\frac{3R}{4}$. We select a cut-off function  $\eta$
$\in C^2(\overline{B}(x_0, R'))$, $0 < \eta < 1$, with $\eta(x_0)=1$ and $\eta = 0$ for $|x-x_0| = R'$,  such that

\begin{eqnarray}
\left.
\begin{array}{rr}
|\nabla \eta|\leq C R^{-1} \eta ^{\alpha}\\
|D^2 \eta|+\eta^{-1}|\nabla \eta|^2\leq C R^{-2} \eta^{\alpha}
\end{array}
\right\}
\quad\text{for}\quad |x-x_0|<R'
\end{eqnarray}
with $C = C (\alpha) > 0$  (see \cite{Zhan} for an example of such function).

First let us state the following lemma.

\begin{lem}\label{tech}
Let $u_0$, $u$ be as in Theorem \ref{profil}.  We denote $w= |\nabla u|^2$ and $z=\eta w$.
Then  at any point $(x_1, t_1)\in Q^{t_0}_{T, R'}$ such that  $|\nabla u(x_1, t_1)|>0$,  $z$ is smooth and  satisfies the following differential inequality
\begin{eqnarray*}
\mathcal{L} z+ C z^{\frac{2q-p+2}{2}}\leq C\left(\dfrac{\left\|u_0\right\|_{\infty}}{t_0}\right)^{\frac{2q-p+2}{q}} + CR^{-\frac{2q-p+2}{q-p+1}},
\end{eqnarray*}
where
\begin{eqnarray}
\mathcal{L} z=\partial_t z- \mathcal{A} z-H \cdot\nabla z,\\
\mathcal{A} z= |\nabla u|^{p-2} \Delta z + (p-2) |\nabla u|^{p-4}(\nabla u)^t D^2 z \nabla u,
\end{eqnarray}
 $H$ is defined by (\ref{apender}) and $C=C(p, q, N)>0$.
 \end{lem}

\textbf{Proof of lemma \ref{tech}}
We know that  a solution  $u$ of \eqref{equi} is smooth at points where $|\nabla u|>0$ \cite{bartier}.
More precisely, we know that  $\nabla u\in C^{2,1}$ in a neighborhood of  such points and hence we can differentiate the equation. As observed in \cite{bartier},  $w=|\nabla u|^2$   satisfies  the following differential equation:
\begin{eqnarray*}
 \partial_t w-\mathcal{A}w =-2 |\nabla u|^{p-2}|D^2 u|^2 +H\cdot\nabla w
\end{eqnarray*}

Indeed, for $i=1,\cdots, N$, put  $u_i=\dfrac{\partial u}{\partial x_i}$ and $w_i=\dfrac{\partial w}{\partial x_i}$.    Differentiating \eqref{equi} in $x_i$, we have
\begin{eqnarray}\label{julien}
\partial_t u_i-|\nabla u|^{p-2}\Delta u_i-\dfrac{p-2}{2} |\nabla u|^{p-4} \sum_{j=1}^N\dfrac{\partial w_i}{\partial x_j} u_j-\dfrac{p-2}{2} |\nabla u|^{p-4} \sum_{j=1}^{N} w_j \dfrac{\partial u_i}{\partial x_j} \nonumber\\
\qquad\qquad =\dfrac{q}{2}  w^{\frac{q-2}{2}} w_i+\dfrac{p-2}{2}   w^{\frac{p-4}{2}} w_i\Delta u+\dfrac{(p-2)(p-4)}{4} w^{\frac{p-6}{2}} (\nabla u\cdot\nabla w ) w_i.
\end{eqnarray}
Multiplying \eqref{julien} by $2u_i$, summing up, and using $\Delta w=2\nabla u \cdot\nabla (\Delta u)+ 2 |D^2 u|^2$, we deduce that
\begin{equation}
\mathcal{L} w=-2 w^{\frac{p-2}{2}}|D^2 u|^2,
\end{equation}
where
\begin{align}
& \mathcal{L}w: = \partial_t w- |\nabla  u|^{p-2} \Delta w -(p-2)|\nabla u|^{p-4} (\nabla u)^t D^2 w \nabla u- H \cdot\nabla w,\\
& H:= \left[ (p-2) w^{\frac{p-4}{2}}\Delta u+\dfrac{(p-2)(p-4)}{2} w^{\frac{p-6}{2}} \nabla u \cdot\nabla w+q w^{\frac{q-2}{2}}\right]\nabla u \nonumber \\
 &\quad\quad+\dfrac{p-2}{2} w^{\frac{p-4}{2}} \nabla w. \label{apender}
\end{align}
Setting $z=\eta w$, we get
\begin{eqnarray*}
 \mathcal{L}z&=&\eta \mathcal{L} w+ w \mathcal{L}\eta-2 w^{\frac{p-2}{2}}\nabla\eta\cdot\nabla w
 -2(p-2) w^{\frac{p-4}{2}} ( \nabla \eta \cdot\nabla u)(\nabla w\cdot\nabla u). 
\end{eqnarray*}

Now we shall estimate the different terms. In what follows $\delta_i>0$ can be chosen arbitrarily small.

\begin{itemize}
\item  Estimate of $|2 w^{\frac{p-2}{2}}\nabla\eta\cdot\nabla w|$.

Using Young's inequality, we have 

\begin{equation}
|2 w^{\frac{p-2}{2}}\nabla\eta\cdot\nabla w|\leq w^{\frac{p-2}{2}} \left[C\eta^{-1}|\nabla\eta|^2w+\delta_1 \eta |D^2 u|^2\right],
\end{equation}
where we used the fact that $\nabla w=2 D^2 u \nabla u$.

\item Estimate of  $|2 (p-2)w^{\frac{p-4}{2}} (\nabla \eta\cdot \nabla u) (\nabla w\cdot\nabla u)|$.
\begin{equation}\label{previous}
 |2(p-2)w^{\frac{p-4}{2}}( \nabla \eta\cdot\nabla u)( \nabla w \cdot\nabla u)|
\leq  w^{\frac{p-2}{2}}\Big[ C\eta^{-1}|\nabla\eta|^2 w+\delta_2 \eta |D^2 u|^2\Big].
\end{equation}
\item Estimate of $|w \, H\cdot\nabla \eta|$.
\begin{align}
 |w \, H \cdot\nabla\eta|&\leq
\underbrace{w^{\frac{p-2}{2}}\left(C\eta^{-1}|\nabla\eta|^2 w+\delta_3 [D^2u|^2\eta\right)}_{(1)} +
\underbrace{w^{\frac{p-2}{2}}\left(C\eta^{-1}|\nabla\eta|^2 w+\delta_4 [D^2u|^2\eta\right)}_{(2)} \nonumber \\
&+\underbrace{w^{\frac{p-2}{2}}\left(C\eta^{-1}|\nabla\eta|^2 w+\delta_5 [D^2u|^2\eta\right)}_{(3)} +C w^{\frac{q+1}{2}}|\nabla \eta|.
\end{align} 

(1) comes from an estimate based on Young's inequality of $ w^{\frac{p-2}{2}}\Delta u(\nabla  u\cdot\nabla\eta)$,
(2) comes from \eqref{previous} and (3) comes from an estimate of $ w^{\frac{p-2}{2}}\nabla w\cdot\nabla\eta$.

\end{itemize}
Finally choosing  $\delta_i$ such that $-2+\delta_1+\delta_2+\delta_3+\delta_4+\delta_5=-1$, we arrive at 
\begin{eqnarray*}
\mathcal{L}z+\eta w^{\frac{p-2}{2}} |D^2 u|^2\leq& C(p,q, N) w^{\frac{p}{2}} \big[|D^2 \eta|+|\Delta \eta|+\eta^{-1}|\nabla\eta|^2\big]+|\nabla\eta| w^{\frac{q+1}{2}}.
\end{eqnarray*}
Using the properties of the cut-off function $\eta$, we get
\begin{eqnarray}\label{jon}
\mathcal{L} z +\eta w^{\frac{p-2}{2}} |D^2 u|^2&\leq& C(p,q,N)  R^{-2} \eta^{\alpha}w^{\frac{p}{2}} +C(p,q,N) R^{-1} \eta^{\alpha}  w^{\frac{q+1}{2}}. 
\end{eqnarray}
Using the result of Theorem \ref{rek}, we shall estimate $ |\nabla u|^{p-2}|D^2 u|^2$ in terms of a power of $w$.
For $(x_1,t_1)\in Q^{t_0}_{ T, R'}$ such that $|\nabla u(x_1, t_1)|>0$, we have 
\begin{eqnarray*}
|\nabla u(x_1, t_1)|^{q}&=& \partial_t u(x_1, t_1)-\text{div} \left(|\nabla u|^{p-2}\nabla u (x_1, t_1)\right)\\
&\leq & \dfrac{\left\|u_0\right\|_{\infty}}{(p-2)t_0} + (p-2+\sqrt{N}) |\nabla u|^{p-2}|D^2 u (x_1,t_1)|.
\end{eqnarray*}
Hence
$$\dfrac{1}{2 (p-2+\sqrt{N})^2} |\nabla u(x_1, t_1)|^{2q}\leq \left(\dfrac{\left\|u_0\right\|_{\infty}}{(p-2)(p-2+\sqrt{N})t_0}\right)^2 +|\nabla u|^{2p-4}|D^2 u(x_1, t_1)|^2.$$
There are two cases:
\begin{equation*}
\begin{array}{ll}
 \text{either} &\dfrac{1}{2 (p-2+\sqrt{N})^2} |\nabla u(x_1, t_1)|^{2q}\leq 2\left(\dfrac{\left\|u_0\right\|_{\infty}}{(p-2)(p-2+\sqrt{N})t_0}\right)^2,\\
 \text{or} &\dfrac{1}{2 (p-2+\sqrt{N})^2} |\nabla u(x_1, t_1)|^{2q-p+2}\leq 2 |\nabla u|^{p-2}|D^2 u(x_1, t_1)|^2.
\end{array}
\end{equation*}
In both cases we arrive at
$$\frac{1}{C(N,p)} |\nabla u(x_1, t_1)|^{2q-p+2}\leq C(p, q, N) \left(\dfrac{\left\|u_0\right\|_{\infty}}{t_0}\right)^{\frac{2q-p+2}{q}}+|\nabla u|^{p-2}|D^2 u(x_1, t_1)|^2.$$
Using this inequality, it follows from  \eqref{jon} that, at $(x_1, t_1)$,

\begin{align*}
\mathcal{L} z+\frac{1}{C(N,p)}\eta |\nabla u|^{2q-p+2}&\leq C(p, q, N) \left(\dfrac{\left\|u_0\right\|_{\infty}}{t_0}\right)^{\frac{2q-p+2}{q}} + C R^{-2}\eta^{\alpha} w^{\frac{p}{2}}+ C R^{-1}\eta^{\alpha} w^{\frac{q+1}{2}}.
\end{align*}

We take  $\alpha=\dfrac{q+1}{2q-p+2}\in(0,1)$ (since $q>p-1)$.
Using Young's inequality and $\eta\leq 1$, we get
$$\mathcal{L} z+\frac{1}{ C(N,p)} \eta |\nabla u|^{2q-p+2}\leq C(p, q, N) \left(\dfrac{\left\|u_0\right\|_{\infty}}{t_0}\right)^{\frac{2q-p+2}{q}} + CR^{-\frac{2q-p+2}{q-p+1}}+ \frac{1}{2 C(N,p)}  \eta |\nabla u|^{2q-p+2} .$$
Hence
\begin{equation}\label{ouf}
\mathcal{L} z+\frac{1}{2 C(N,p)} z^{\frac{2q-p+2}{2}}\leq C(p, q, N) \left(\dfrac{\left\|u_0\right\|_{\infty}}{t_0}\right)^{\frac{2q-p+2}{q}}  + CR^{-\frac{2q-p+2}{q-p+1}}.
\end{equation}

\bigskip

\subsection*{Proof of theorem \ref{profil}}

First let us note that  by the proof of the local existence there exists $t_0\in \Big(0, T_{max}(u_0)\Big)$  with $t_0= t_0(M, p, q, N, \left\|g\right\|_{C^2})$, such that 
\begin{equation}\label{utile}
\underset{0\leq t\leq t_0}{\text{sup}} \left\|\nabla u\right\|_{L^{\infty}}\leq C( p, q,  \Omega, M, \left\|g\right\|_{C^2}).
\end{equation}
 We also know that $\nabla u$ is a locally $\mathrm{H\ddot{o}lder}$ continuous function and thus  $z$ is a continuous function on $\overline{B(x_0, R')}\times [t_0, T]=\overline{Q}$, for  any $T< T_{max(u_0)}$. Therefore, unless $z\equiv 0$ in $\overline{Q}$,  $z$  must reach a positive maximum at some point $(x_1,t_1)\in \overline{B(x_0, R')}\times [t_0, T]$.  
Since $z=0$ on $\partial B_{R'}\times [t_0, T]$, we deduce that $x_1\in B_{R'}$. Therefore $\nabla z(x_1, t_1)=0$ and $D^2 z(x_1,t_1)\leq 0$.  Now we have either $t_1=t_0$, or $t_0< t_1\leq T$.
If $t_1=t_0$, then 
$$z(x_1, t_1)\leq \left\|\nabla u(t_0)\right\|_{L^{\infty}}^2\leq C( p, q,  \Omega, M, \left\|g\right\|_{C^2}).$$
If $t_0< t_1\leq T$, we have $\partial_t z(x_1,t_1)\geq 0$ and therefore $\mathcal{L} z\geq0$.
Using \eqref{ouf} we arrive at
\begin{equation}\label{ouff}
\frac{1}{2 C(N,p)} z(x_1,t_1)^{\frac{2q-p+2}{2}}\leq  C(p, q, N) \left(\dfrac{\left\|u_0\right\|_{\infty}}{t_0}\right)^{\frac{2q-p+2}{q}}+ CR^{-\frac{2q-p+2}{q-p+1}},
\end{equation}
 that is 
\begin{equation}\label{oug}
 \sqrt{ z(x_1,t_1)}\leq C(p, q, N) \left(\dfrac{\left\|u_0\right\|_{\infty}}{t_0}\right)^{\frac{1}{q}} + C(p,q, N)R^{-\frac{1}{q-p+1}}.
\end{equation}
Since $z(x_0,t)\leq z(x_1, t_1)$ and $\eta(x_0)=1$, we get
 $$|\nabla u(x_0,t)|\leq C(p, q, N) \left(\dfrac{\left\|u_0\right\|_{\infty}}{t_0}\right)^{\frac{1}{q}} + C(p,q, N)R^{-\frac{1}{q-p+1}}\quad\text{for}\quad t\in  [t_0, T].$$
The proof of \eqref{profil} follows by  taking $R=\delta(x_0)$, letting $T\to T_{max}(u_0)$ and using \eqref{utile}.

\section{Blow-up criterion: proof of Theorem \ref{blo}}
 
Assume that $T_{max}(u_0)=\infty$, taking $\varphi_1^{\alpha}$ as  test-function, we have for any $\tau>0$
\begin{equation}\label{formule}
\int_0^{\tau}  \int_{\Omega} u_t\, \varphi_1^{\alpha}\, dx dt =\int_0^{\tau}\int_{\Omega} |\nabla u|^q\,\varphi_1^{\alpha} \,dx dt -\alpha \int_0^{\tau} \int_{\Omega} |\nabla u|^{p-2} \,\varphi_1^{\alpha-1}\, \nabla u\cdot\nabla \varphi_1 \,dx dt.
\end{equation}

Set $y(t) =\int_{\Omega} u(t)\,\varphi_1^{\alpha} \,dx$. Since by  definition $u_t\in L^2_{loc}((0,\infty); L^2(\Omega))$, we have $y\in  W^{1,1}_{loc}(0,\infty)$ and $y'(t)=\int_{\Omega} u_t\,\varphi_1^{\alpha}\, dx$. Differentiating \eqref{formule} with respect to $\tau$ we have, for a.e. $\tau>0$
\begin{equation}
y'(\tau)=\int_{\Omega} |\nabla u(\tau)|^q\,\varphi_1^{\alpha} \,dx -\alpha  \int_{\Omega} |\nabla u(\tau)|^{p-2} \,\varphi_1^{\alpha-1}\, \nabla u(\tau)\cdot\nabla \varphi_1 \,dx.
\end{equation}
Assume that \textbf{$\alpha >\frac{p-1}{(q-p+1)}$}. Since $q> p>1$ and $\left\| \nabla \varphi_1\right\|_{\infty}\leq C'$, using $\mathrm{H\ddot{o}lder}$  and Young inequalities we get:
 \begin{eqnarray*}
\alpha \int_{\Omega}|\nabla  u(\tau)|^{p-2}\, \varphi_1^{\alpha-1}\,\nabla u(\tau)\cdot\nabla \varphi_1\, dx &\leq &\frac{1}{2} \int_{\Omega}|\nabla u(\tau)|^q\, \varphi_1^{\alpha}\, dx  +C \int_{\Omega} \varphi_1^{\alpha- q/(q-p+1)} \,dx\\
 &\leq&\frac{1}{2} \int_{\Omega}|\nabla u(\tau)|^q \varphi_1^{\alpha} dx  +C.
\end{eqnarray*}
Here we used the fact that  $\int_{\Omega} \varphi_1^{-l} dx <\infty $ for $l<1$. 
Therefore  $$y'(\tau) \geq \frac{1}{2} \int_{\Omega}|\nabla u(\tau)|^q\varphi_1^{\alpha} dx -C.$$
Assuming that $\alpha < q-1$, we get
\begin{align*}
\int_{\Omega} |\nabla u(\tau)|\, dx=\int_{\Omega} |\nabla u(\tau)|\, \varphi_1^{\frac{\alpha}{q}}\,\varphi_1^{-\frac{\alpha}{q}}\,dx&\leq\left( \int_{\Omega} |\nabla u(\tau)|^q\, \varphi_1^{\alpha}\,dx \right)^{1/q} \left( \int_{\Omega} \varphi_1^{\frac{-\alpha}{q-1}}\,dx\right)^{\frac{q-1}{q}}\\
&\leq C \left( \int_{\Omega} |\nabla u(\tau)|^q\, \varphi_1^{\alpha}\,dx \right)^{1/q}.
\end{align*}
On the other hand using that $ \int_{\Omega} |u(\tau)| \,dx\leq C \left\|u\right\|_{L^{\infty}(\partial\Omega)}+ C\int_{\Omega} |\nabla u(\tau)| \,dx$,  we have
$$\int_{\Omega} u(\tau)\, \varphi_1^{\alpha}\,dx\leq \left\|\varphi_1^{\alpha}\right\|_{\infty}\, \int_{\Omega} u(\tau)\,dx\leq  C+C\,\int_{\Omega} |\nabla u(\tau)|\,dx .$$
Combining these two inequalities we arrive at
$$\int_{\Omega} |\nabla u(\tau)|^q\, \varphi_1^{\alpha}\,dx \geq  C \left( \int_{\Omega} u(\tau)\, \varphi_1^{\alpha}\,dx \right)^{q}- C.$$
Finally we get the blow-up inequality
$$y'(\tau)\geq C_1 y(\tau)^q - C_2, \quad \text{ for a.e.}\quad \tau>0,$$
with $C_1=C_1(p, q, \Omega)>0$ and $C_2=C_2(p,q,\alpha, \Omega, \left\|g\right\|_{\infty})$.

\begin{rem}\label{secblo}
Instead of assuming that  $\int_{\Omega}u_0 \phi_1^{\alpha} dx$ is large in Theorem \ref{blo}, it would be sufficient to assume that $\left\|u_0\right\|_r$ is large for some $r\in [1,\infty)$. In fact,
assuming without loss of generality $r \geq (2q-p)/(q-p)$ and denoting $y(t) =\int_{\Omega}u^r(t) dx$, the Poincar$\grave{e}$ and H$\ddot{o}$lder inequalities can be used in order to prove the blow-up inequality $y'\geq C_1 y^{(q+r-1)/r}-C_2$ (see \cite{oups})\rmfamily.
\end{rem}

\textbf{Acknowledgements.} The author would like to thank Professor Ph. Souplet
for useful suggestions during the preparation of this paper.

\nocite{*}
\bibliography{bibliomai}
\bibliographystyle{amsplain}

\end{document}

%% file: mondessin.pgf
\begin{pgfpicture}{-3.5556cm}{-3.5556cm}{4.4444cm}{3.5556cm}%
\pgfsetxvec{\pgfxy(1.7778,0)}
\pgfsetyvec{\pgfxy(0,1.7778)}
\pgfsetstrokecolor{rgb,1:red,0;green,0;blue,0}
\pgfsetlinewidth{0.4pt} 
\pgfsetroundjoin \pgfpathmoveto{\pgfxy(0,0)}\pgfpathlineto{\pgfxy(0.0095,0.1259)}\pgfpathlineto{\pgfxy(0.0348,0.2532)}
\pgfpathlineto{\pgfxy(0.0711,0.3809)}\pgfpathlineto{\pgfxy(0.1136,0.5085)}\pgfpathlineto{\pgfxy(0.1575,0.6349)}
\pgfpathlineto{\pgfxy(0.1981,0.7595)}\pgfpathlineto{\pgfxy(0.2305,0.8815)}\pgfpathlineto{\pgfxy(0.2377,0.9253)}
\pgfpathlineto{\pgfxy(0.1628,0.8597)}\pgfpathlineto{\pgfxy(0.0762,0.7659)}\pgfpathlineto{\pgfxy(0,0.6634)}
\pgfpathlineto{\pgfxy(-0.0648,0.5533)}\pgfpathlineto{\pgfxy(-0.1176,0.437)}\pgfpathlineto{\pgfxy(-0.1577,0.3157)}
\pgfpathlineto{\pgfxy(-0.1847,0.1909)}\pgfpathlineto{\pgfxy(-0.1983,0.0639)}\pgfpathlineto{\pgfxy(-0.1983,-0.0639)}
\pgfpathlineto{\pgfxy(-0.1847,-0.1909)}\pgfpathlineto{\pgfxy(-0.1577,-0.3157)}\pgfpathlineto{\pgfxy(-0.1176,-0.437)}
\pgfpathlineto{\pgfxy(-0.0648,-0.5533)}\pgfpathlineto{\pgfxy(0,-0.6634)}\pgfpathlineto{\pgfxy(0.0762,-0.7659)}
\pgfpathlineto{\pgfxy(0.1628,-0.8597)}\pgfpathlineto{\pgfxy(0.2377,-0.9253)}\pgfpathlineto{\pgfxy(0.2305,-0.8815)}
\pgfpathlineto{\pgfxy(0.1981,-0.7595)}\pgfpathlineto{\pgfxy(0.1575,-0.6349)}\pgfpathlineto{\pgfxy(0.1136,-0.5085)}
\pgfpathlineto{\pgfxy(0.0711,-0.3809)}\pgfpathlineto{\pgfxy(0.0348,-0.2532)}\pgfpathlineto{\pgfxy(0.0095,-0.1259)}
\pgfpathlineto{\pgfxy(0,0)}
\pgfsetfillcolor{rgb,1:red,0.502;green,0.502;blue,0.502}
\pgffill
\pgfsetlinewidth{2pt} 
\pgfpathmoveto{\pgfxy(0,0)}
\pgfpathcurveto{\pgfxy(0,0.3333)}{\pgfxy(0.2195,0.6892)}{\pgfxy(0.25,1)}
\pgfpathcurveto{\pgfxy(0.2805,1.3108)}{\pgfxy(0.122,1.5765)}{\pgfxy(-0.25,1.6)}
\pgfpathcurveto{\pgfxy(-0.622,1.6235)}{\pgfxy(-1.2073,1.4047)}{\pgfxy(-1.5,0.9)}
\pgfpathcurveto{\pgfxy(-1.7927,0.3953)}{\pgfxy(-1.7927,-0.3953)}{\pgfxy(-1.5,-0.9)}
\pgfpathcurveto{\pgfxy(-1.2073,-1.4047)}{\pgfxy(-0.622,-1.6235)}{\pgfxy(-0.25,-1.6)}
\pgfpathcurveto{\pgfxy(0.122,-1.5765)}{\pgfxy(0.2805,-1.3108)}{\pgfxy(0.25,-1)}
\pgfpathcurveto{\pgfxy(0.2195,-0.6892)}{\pgfxy(0,-0.3333)}{\pgfxy(0,0)}
\pgfstroke
\pgfsetlinewidth{0.6pt} 
\pgfellipse[stroke]{\pgfxy(1,0)}{\pgfxy(0,1)}{\pgfxy(-1,0)}
\pgfsetlinewidth{1pt} 
\pgfellipse[stroke]{\pgfxy(1,0)}{\pgfxy(0,1.2)}{\pgfxy(-1.2,0)}
\pgfputat{\pgfxy(1,-0.1406)}{\pgftext[top]{\color{rgb,1:red,0;green,0;blue,0}\small $x_1$}}\pgfstroke
\pgfputat{\pgfxy(-1.3594,0.9)}{\pgftext[left]{\color{rgb,1:red,0;green,0;blue,0}\small $\partial \Omega$}}\pgfstroke
\pgfputat{\pgfxy(-0.5398,0.3398)}{\pgftext[right,bottom]{\color{rgb,1:red,0;green,0;blue,0}\small $\Gamma$}}\pgfstroke
\pgfsetlinewidth{0.4pt} 
\pgfsetfillcolor{rgb,1:red,0;green,0;blue,0}\pgfellipse[fillstroke]{\pgfxy(0,0)}{\pgfxy(0,0.0277)}{\pgfxy(-0.0277,0)}
\pgfputat{\pgfxy(0.0994,-0.0994)}{\pgftext[left,top]{\color{rgb,1:red,0;green,0;blue,0}\small $\widetilde{x_0}$}}\pgfstroke
\pgfellipse[fillstroke]{\pgfxy(-0.12,0)}{\pgfxy(0,0.0277)}{\pgfxy(-0.0277,0)}
\pgfputat{\pgfxy(-0.12,0.1406)}{\pgftext[bottom]{\color{rgb,1:red,0;green,0;blue,0}\small $x_0$}}\pgfstroke
\pgfsetarrows{to-to}
\pgfxyline(1,0)(0,0)
\pgfsetarrows{-}
\pgfputat{\pgfxy(0.5,0.1406)}{\pgftext{\color{rgb,1:red,0;green,0;blue,0}\small $\rho$}}\pgfstroke
\pgfsetarrows{to-to}
\pgfxyline(1,1)(1,1.2)
\pgfsetarrows{-}
\pgfputat{\pgfxy(0.8594,1.1)}{\pgftext{\color{rgb,1:red,0;green,0;blue,0}\small $\eta$}}\pgfstroke
\pgfsetarrows{to-to}
\pgfxyline(1,0)(0.017,-0.6883)
\pgfsetarrows{-}
\pgfputat{\pgfxy(0.5892,-0.4593)}{\pgftext{\color{rgb,1:red,0;green,0;blue,0}\small $\eta+\rho$}}\pgfstroke
\pgfsetarrows{-to}
\pgfxyline(-0.5,0.3)(-0.08,-0.2)
\pgfsetarrows{-}
\end{pgfpicture}%